\documentclass{scrartcl}

\KOMAoptions{paper=a4}
\KOMAoptions{fontsize=12pt}
\KOMAoptions{DIV=calc} 

\usepackage[utf8]{inputenc} 
\usepackage[T1]{fontenc} 
\usepackage[ngerman,english]{babel} 
\usepackage{csquotes} 
\usepackage{hyphenat} 

\usepackage{kpfonts}

\usepackage{cabin}

\usepackage{microtype} 

\addtokomafont{title}{\rmfamily}

\usepackage{enumitem}
\setlist[enumerate]{label*=(\alph*),ref=(\alph*),itemsep=0pt,topsep=5pt}
\setlist[itemize]{itemsep=0pt,topsep=5pt}

\KOMAoptions{DIV=calc}

\usepackage{booktabs} 

\usepackage{graphicx} 
\graphicspath{{pic/}} 

\usepackage[format=plain,labelfont={sf,footnotesize},textfont=small,margin=12pt]{caption} 
\DeclareCaptionLabelSeparator{MySpace}{\enskip } 
\usepackage[style=alphabetic,
						sorting=nyt,
						maxnames=4,
						backend=biber,
						doi=true,
						isbn=false,
						clearlang=true]
						{biblatex}

\bibliography{C:/Users/j.rau/Dropbox/Bibliography/biblio}

\DeclareFieldFormat*{title}{\textit{#1}} 
\renewbibmacro{in:}{} 
\DeclareFieldFormat{journaltitle}{\iffieldequalstr{journaltitle}{ArXiv e-prints}{Preprint}{#1}} 
\DeclareFieldFormat{booktitle}{#1} 
\DeclareRedundantLanguages{english,English}{english,german,ngerman,french}
\DeclareSourcemap{ 
  \maps[datatype=bibtex]{
    \map{
      \step[fieldset=addendum,null]
      \step[fieldset=eprintclass,null]
    }  
  }
}

\usepackage{hyperref} 
\usepackage{xcolor} 
\hypersetup{breaklinks=true} 
\definecolor{darkblue}{RGB}{0,0,170}
\definecolor{darkred}{RGB}{200,0,0}
\hypersetup{citecolor=darkblue, urlcolor=darkblue, linkcolor=darkblue, colorlinks=true} 
\usepackage[all]{hypcap} 

\usepackage{aliascnt} 
\usepackage[amsmath,thmmarks]{ntheorem} 

\newtheorem {theorem}{Theorem}[section]

\newaliascnt{proposition}{theorem}
\newtheorem {proposition}[proposition]{Proposition}
\aliascntresetthe{proposition}

\newaliascnt{lemma}{theorem}
\newtheorem {lemma}[lemma]{Lemma}
\aliascntresetthe{lemma}

\newaliascnt{corollary}{theorem}

\aliascntresetthe{corollary}

\newaliascnt{conjecture}{theorem}

\aliascntresetthe{conjecture}

\theorembodyfont{\upshape}

\newaliascnt{definition}{theorem}
\newtheorem {definition}[definition]{Definition}
\aliascntresetthe{definition}

\newaliascnt{example}{theorem}

\aliascntresetthe{example}

\newaliascnt{exercise}{theorem}

\aliascntresetthe{exercise}

\newaliascnt{goal}{theorem}

\aliascntresetthe{goal}

\newaliascnt{construction}{theorem}

\aliascntresetthe{construction}

\theoremheaderfont{\itshape}
\theorembodyfont{\upshape}

\newaliascnt{remark}{theorem}
\newtheorem {remark}[remark]{Remark}
\aliascntresetthe{remark}

\newaliascnt{convention}{theorem}

\aliascntresetthe{convention}

\newaliascnt{notation}{theorem}

\aliascntresetthe{notation}

\theoremstyle {nonumberplain}
\theoremheaderfont{\itshape}
\theorembodyfont{\upshape}
\theoremsymbol{\ensuremath{_\blacksquare}}

\newtheorem {proof}{Proof}

\usepackage {amsmath} 

\newcommand{\F}{{\mathbf F}}

\newcommand{\R}{{\mathbf R}}

\newcommand{\Z}{{\mathbf Z}}

\newcommand{\CP}{\mathbf{CP}}

\newcommand{\CC}{{\mathcal C}}
\newcommand{\FF}{{\mathcal F}}
\newcommand{\GG}{{\mathcal G}}

\DeclareMathOperator{\Tr}{Tr}
\DeclareMathOperator{\BM}{BM}
\DeclareMathOperator{\rk}{rk}

\DeclareMathOperator{\Star}{Star}

\DeclareMathOperator{\Aut}{Aut}

\DeclareMathOperator{\Fix}{Fix}

\DeclareMathOperator{\gap}{gap}

\DeclareMathOperator{\GL}{GL}
\DeclareMathOperator{\cl}{cl}

\begin {document}

\title {On the tropical \\ Lefschetz-Hopf trace formula}
\author {Johannes Rau%
  \thanks{
	  The author is supported by the FAPA project \enquote{Matroids in tropical geometry} 
		at Universidad de los Andes, Colombia.
		\emph{Author e-mail address}: \href{mailto:j.rau@uniandes.edu.co}{j.rau@uniandes.edu.co}
	}
}
\date{}

\maketitle

\begin{abstract}
	\noindent
  In this follow-up to \cite{Rau-TropicalPoincareHopf}, our main result 
	is a tropical Lefschetz-Hopf trace formula for 
	matroidal automorphisms. We show that both sides 
	of the formula are equal to 
	the (generalized)
	beta invariant of the lattice of fixed flats. 
\end{abstract}


\section{Introduction}

In \cite{Rau-TropicalPoincareHopf}, the author proves a tropical version
of the Poincaré-Hopf theorem and conjectures a tropical version
of a tropical Lefschetz-Hopf trace formula.
The main result of this paper is the proof of the tropical trace formula
in the case of \emph{matroidal} automorphisms. In this case, the formula can
be refined by giving a third description of the value
in terms of the (generalized) beta invariant 
of the lattice of fixed flats. 

\begin{theorem} \label{thm:traceformulamatroidauts}
  Let $M$ be a loopless matroid and let $\Psi \colon \Sigma_M \to \Sigma_M$
	be a matroidal automorphism. Then
	\begin{equation} \label{eq:traceformula} 
		\deg (\Gamma_{\Psi} \cdot \Delta) = (-1)^n \beta(\Fix(L(M))) = 
		  \sum_{p=0}^n (-1)^{p} \Tr(\Psi_*, \F_p(\Sigma_M)).
	\end{equation}
\end{theorem}

The theorem is an instance of the conjecture
stated in \cite[Conjecture 1.2]{Rau-TropicalPoincareHopf}
using standard instead of Borel-Moore tropical homology
(cf.\ \cite[Remark 1.8]{Rau-TropicalPoincareHopf}),
since $H_{p,0}(\Sigma_M) = \F_p(\Sigma_M)$ and
$H_{p,q}(\Sigma_M) = 0$ for $q >0$. We actually
hope that the theorem can be used as the local ingredient
to prove more general cases of 
\cite[Conjecture 1.2]{Rau-TropicalPoincareHopf}.

Let us quickly review the ingredients of the theorem 
(precise definitions follow in the later sections).
In the following, $M$ is a loopless matroid of rank $n+1$ on 
the set $E = \{0, \dots, N\}$.
We denote by $\Sigma_M \subset \R^N$ the 
(projective) matroid fan of $M$. 
A \emph{matroidal automorphism} $\Psi \colon \Sigma_M \to \Sigma_M$
is a tropical automorphism of $\Sigma_M$ which is induced by a
matroid automorphism $\psi \colon M \to M$. More precisely, the relationship 
is given by
\[
  \Psi(v_S) = v_{\psi(S)}
\]
for any $S \subset E$ and associated indicator vector $v_S$. 

Given $\Psi$, we denote by $\Gamma_\Psi, \Delta \subset \Sigma_M \times \Sigma_M$
the graph of $\Psi$ and the diagonal of $\Sigma_M$, respectively. They carry the structure 
of tropical cycles. Using the tropical intersection theory for matroid fans 
constructed in \cite{Sha-TropicalIntersectionProduct, FR-DiagonalTropicalMatroid},
we can define their intersection product $\Gamma_\Psi \cdot \Delta$.  
The degree of this product is the intersection-theoretic side of the tropical
trace formula.

For $p = 0, \dots, n$, the framing group $\F_p(\Sigma_M) \subset \bigwedge^p \R^N$ 
is the vector space generated by wedges of $p$ vectors contained in the same 
cone of $\Sigma_M$ \cite{MZ-TropicalEigenwaveIntermediate, IKMZ-TropicalHomology}. 
Since $\Psi$ is the restriction of a linear map on $\R^N$, 
$\Psi$ induces a map $\Psi_* \colon \bigwedge^p \R^N \to \bigwedge^p \R^N$
which can be restricted to $\Psi_* \colon \F_p(\Sigma_M) \to \F_p(\Sigma_M)$.
We denote the trace of this automorphism by $\Tr(\Psi_*, \F_p(\Sigma_M))$. 
The alternating sum of traces is the trace side of the tropical trace formula.

To connect the two sides, we use an intermediate, combinatorially defined
invariant which is the (generalized) \emph{beta invariant} of a lattice
with rank function. In our case, the lattice in question is
the lattice of flats $F \in L(M)$ which are fixed by $\psi$,
\[
  \Fix(L(M)) = \{F \in L(M) : F = \psi(F)\}
\]
The beta invariant is defined using the Möbius function $\mu^\psi$ of 
$\Fix(L(M))$ and the rank function $\rk$ of $M$ (or $L(M)$) by
\[
  \beta(\Fix(L(M))) = 
	  (-1)^{n+1} \sum_{\substack{F \in L(M) \\ F = \psi(F)}} \mu^\psi(\emptyset, F) \rk(F).
\]
To prove the trace formula, we will show that both sides agree with the beta invariant
(up to sign). Note that this is a generalization of the situation
in \cite{Rau-TropicalPoincareHopf}, where the two sides of the Poincaré-Hopf
theorem are shown to be equal to the (ordinary) beta invariant of $M$. 
Even though we only use elementary properties of the
generalized beta invariant and the fixed lattice $\Fix(L(M))$,
we hope that this results encourages the further study 
of their combinatorial properties.

\begin{remark} 
	Let us emphasize that $\Sigma_M$ denotes the projective matroid fan 
	in $\R^N \cong \R^{N+1}/\R \mathbf{1}$. 
	The reasons for preferring it over the affine matroid fan $\Sigma'_M$ are essentially the ones
	mentioned in the beginning of \cite[Section 2.3]{Rau-TropicalPoincareHopf}.
	First, the analogous intersection product 
	$\Gamma_{\Psi'} \cdot \Delta' \subset \Sigma'_M \times \Sigma'_M$
	is zero for trivial reasons. Second, when trying to compute $\Gamma_{\Psi} \cdot \Delta$
	by taking pull-backs to $\Sigma'_M \times \Sigma'_M$, we are leaving the class of matroid fans 
	which makes certain combinatorial arguments less transparent and natural.
	To overcome this, one would have to perform a break of symmetry analogous 
	to the one we will encounter for example in \autoref{eq:functiongidehom}, but without
	benefits. 
\end{remark}

\begin{remark} 
  In general, the notion of (auto-)morphism commonly used in tropical geometry
	is based on locally $\Z$-linear affine maps. 
	This allows for more general automorphisms of matroid fans
	than the matroidal automorphisms considered here. 
	A trivial example is 
	$\R^N = \Sigma_{U_{N+1,N+1}}$ with $\Aut(\R^N) = \GL(N,\Z) \ltimes \R^N$ but
	$\Aut(U_{N+1,N+1}) = S_{N+1}$.
	For this particular case, the trace formula for general automorphisms follows
	from \cite[Section 4.5]{Rau-TropicalPoincareHopf}. 
	For arbitrary matroid fans, the situation is unclear. 
	However, the case of matroidal automorphisms seems to be interesting in its own right,
	since its purely combinatorial nature yields interesting combinatorial 
	features such as the generalized beta-invariants considered here. 
	Moreover, a future goal is to decompose the (general) automorphism group
	of a matroid fan into various subgroups such as the matroidal automorphisms, 
	the ones that fix the coarse subdivision of $\Sigma_M$ but act by $\Z$-linear
	automorphisms on the lineality space, possibly certain exceptional ones, etc., 
	and reduce the general trace formula to considerations for each subgroup. 
	Encouragement for this approach is provided for example by the results in 
	\cite[Theorems 1.2, 6.3, 9.5]{SW-BirationalGeometryMatroids}
	which establish that in certain cases matroidal automorphisms
	together with certain exceptional automorphisms generate
	(a big subgroup of) all automorphisms. 
	So, we hope that the present analysis is not only of interest 
	for its particular combinatorial features, but also can 
	serve in the future to establish more general tropical trace formulas. 
\end{remark}

\paragraph*{Acknowledgements}

I would like to thank Kristin Shaw, Karim Adiprasito and Omid Amini for useful discussions
and feedback on this project. I would like to thank the two anonymous referees for 
careful reading of the manuscript and many helpful comments and corrections.

\section{Preliminaries}

\subsection{Notational summary}

In order to compute the intersection-theoretic side of the trace formula,
we will use the same approach as in \cite{Rau-TropicalPoincareHopf}.
That is to say, we will rely on the expression of the diagonal in 
terms of generic chains of matroids from \cite{FR-DiagonalTropicalMatroid}. 
In order to keep the overlap to a minimum, we will just give a quick summary
of the relevant notation and statements here and refer to reader to
the aforementioned sources for more details.

\bigskip
\begin{tabular}{lp{0.8\linewidth}}
  $E$ & base set $\{0, \dots, N\}$ \\
  $M$ & loopless matroid of rank $n+1$ on $E$ \\
	$\rk$ & rank function of $M$ \\
  $L(M)$ & lattice of flats of $M$ \\
	$\mathbf{1}$ & all one vector $(1,\dots,1) \in \R^E \cong\R^{N+1}$ \\
	$v_S$ & vector in $\R^E$ (and by abuse of notation in $\R^E/\R\mathbf{1}$) 
	whose $i$-th entry is $-1$ if $i \in S$ and $0$ otherwise, for a subset $S \subset E$ \\
	$\FF$ & chain of flats $E = F_0 \supsetneq F_1 \supsetneq \dots \supsetneq F_l \supsetneq F_{l+1} = \emptyset$ (note that we count in decreasing order) \\
	$l(\FF)$ & length of a chain, which is $l$ in the above notation ($E$ and $\emptyset$ are not counted) \\
	$\gap(\FF)$ & gap sequence $(r_0, r_1, \dots, r_l)$ of $\FF$, where $r_i := \rk F_i - \rk F_{i+1}-1$ \\
	$\sigma_\FF$ & cone in $\R^E$ (and $\R^E/\R\mathbf{1}$) generated by the indicator vectors $v_{F_i}$ (and $\R v_E = \R\mathbf{1}$) \\
	$\Sigma'_M$ & \emph{affine} matroid fan sitting in $\R^E$ (the union of cones $\sigma_\FF$ for all chains of flats $\FF$) \\
	$\Sigma_M$ & \emph{projective} matroid fan sitting in $\R^E/\R\mathbf{1}$ \\
	$\F_p(\Sigma_M)$ & framing groups of $\Sigma_M$ (generated by wedges of $p$ vectors contained in the same cone of $\Sigma_M$)
\end{tabular}
\bigskip

Running through all chains $\CC$ of (arbitrary) subsets of $E$, the collection of cones
$\sigma_\CC$ forms a unimodular subdivision of $\R^E$ (and $\R^E/\R\mathbf{1}$) which 
we call the \emph{permutahedral fan} (since it is the normal fan of the permutahedron). 
Throughout the following, matroid fans (as well as all other fans to come) will always 
be represented as subfans of the permutahedral fan. This representation
is called the \emph{fine subdivision} of $\Sigma_M$. 
Since the permutahedral fan is unimodular, in order to describe a piecewise $\Z$-linear function $f$
on it, it suffices to prescribe its values on the indicator vectors $v_S$. We will often use the
shorthand $f(S)$ instead of $f(v_S)$.

Given two matroids $M, Q$ on the ground set $E$, it is obvious that $\Sigma_Q \subset \Sigma_M$
(both as sets and fans) if and only if $L(Q) \subset L(M)$ (or, in matroid terminology, $Q$ is a \emph{quotient} of $M$). 
In such a case, there exists a canonical sequence of matroids $Q = M_0, M_1, \dots, M_s = M$
(called the \emph{generic chain}) 
such that $\rk(M_i) = \rk(Q) + i$ and with rank functions
\begin{equation} \label{eq:intermediatematroids} 
  \rk_{M_i}(S) = \min\{\rk_Q(S) + i, \rk_M(S)\}.
\end{equation}
Moreover, there is a sequence of piecewise $\Z$-linear functions 
(on the permutahedral fan)
$g'_1, \dots, g'_s : \R^{N+1} \to \R$ such that
\begin{equation} \label{eq:intersectfunctiongeneral} 
  \Sigma'_{M_{s-i}} = g'_i \cdot g'_{i-1} \cdots g'_1 \cdot \Sigma'_M.
\end{equation}
If $Q$ and $M$ correspond to the hyperplane arrangements associated
to the projective subspaces $K \subset L \subset \CP^n$,
then the $M_i$ correspond to a generic flag of subspaces $K \subset S_1 \subset \dots \subset S_s = L$
(see also \autoref{lem:hyperplanesection} for a tropical version).
The functions are given by 
\begin{equation} \label{eq:functioningeneral} 
  g'_i(S) = \begin{cases}
	           -1 & \rk_M(S) \geq \rk_Q(S) + s + 1 - i, \\
						  0 & \text{otherwise}.
	         \end{cases}
\end{equation}
Let us note that we use the max-convention here. Moreover, 
given a weighted unimodular fan $\Sigma$ and 
a function $g \colon |\Sigma| \to \R$ which is linear on the cones of $\Sigma$, 
we denote by $g \cdot \Sigma$ the weighted fan defined in
\cite[Definition 3.4]{AR-FirstStepsTropical}.
It is given by the codimension $1$ skeleton of $\Sigma$ equipped 
with the weights
\[
  \omega_{g \cdot \Sigma}(\tau) = \sum_{\tau \subsetneq \sigma \in \Sigma} 
		\omega_\Sigma(\sigma) g(v_{\sigma/\tau}) 
	  - g\left( \sum_{\tau \subsetneq \sigma \in \Sigma} 
			\omega_\Sigma(\sigma) v_{\sigma/\tau} \right).
\] 
Here, $v_{\sigma/\tau}$ denotes the primitive generator
of the unique ray in $\sigma$ not contained in $\tau$ 
(since $\Sigma$ is unimodular). 
See \autoref{rem:weightcomp} for a more detailed analysis of our particular case.

We denote by $M\oplus_0 M$ the parallel connection of $M$ with itself along the element $0$.
Its base set $E \sqcup_0 E$ is the disjoint union of $E$ with itself, but the two zeros identified.
By convention, we write a subset of $E \sqcup_0 E$ as a pair $(F,G), F,G \subset E$ such that
either $0 \in F \cap G$ or $0 \notin F \cup G$. 
We will always identify $\R^E/\R\mathbf{1} = \R^N$ by setting the coordinate corresponding
to $0 \in E$ to zero --- $x_0 = 0$. This induces a natural identification
of 
\[
  \R^{E\sqcup_0 E}/\R\mathbf{1} = \R^{2N} = \R^E/\R\mathbf{1} \times \R^E/\R\mathbf{1}.
\]
Under this identification, we have $\Sigma_{M\oplus_0 M} = \Sigma_M \times \Sigma_M$. 
Using this setup, the diagonal $\Delta \subset \Sigma_M \times \Sigma_M$ can also be represented as a matroid fan.
The rank function of the associated matroid is given by
\[
  \rk_\Delta(F,G) = \rk(F \cup G).
\]

Applying the construction from \autoref{eq:intersectfunctiongeneral} to $\Delta \subset \Sigma_M \times \Sigma_M$,
we obtain (dehomogenised) functions
$g_1, \dots, g_n \colon \R^{2N} \to \R$
given by
\begin{equation} \label{eq:functiongidehom} 
  g_i(F,G) = \begin{cases}
	             -1 & 0 \notin F, \rk(F) + \rk(G) \geq \rk(F \cup G) + n + 1 - i, \\
	             +1 & 0 \in F, \rk(F) + \rk(G) \leq \rk(F \cup G) + n + 1 - i, \\
						    0 & \text{otherwise},
					   \end{cases}
\end{equation}
and such that 
\begin{equation} \label{completeintersection} 
	  \Delta = g_n \cdots g_1 \cdot (\Sigma_M \times \Sigma_M)
\end{equation}
(see \cite[Equation (12) and Proposition 2.6]{Rau-TropicalPoincareHopf}).
Our computation of $\deg(\Gamma_\Psi \cdot \Delta)$ is based on this description of $\Delta$.

\subsection{Hyperplane sections}

A special case of \autoref{eq:intersectfunctiongeneral} that we will use briefly
is when $Q = U_{1,N+1}$ is the unique loopless matroid of rank $1$ with ground set $E$. 
In this case, $\Sigma_Q = \{\mathbf{0}\}$.
Given any loopless matroid $M$, let $M^{\leq i}$ denote the matroid whose 
flats are the flats of $M$ of rank at most $i$ as well as $E$ (called
the $i$-th \emph{truncation} of $M$).
Indeed, this defines a matroid by an obvious check of axioms or by the following observation.

\begin{lemma} \label{lem:genericmatroids}
  Let $M$ be a loopless matroid and set $Q = U_{1,N+1}$.
	Then the intermediate matroids $M_i$ from \autoref{eq:intersectfunctiongeneral}
	are equal to $M^{\leq i}$ for all $i = 0, \dots, s$. 
\end{lemma}

\begin{proof}
  This is straightforward using \autoref{eq:intermediatematroids}
	and $\rk_Q(S) = 1$ for all $S \neq \emptyset$. 
\end{proof}

As mentioned before, the intermediate matroids $M_i$ 
should be thought of as generic hyperplane sections. 
This can be made precise tropically by the following statement,
which is originally proven in \cite[Lemma 5.1]{HK-LogConcavityCharacteristic}.
As a warm-up, we include a short proof based on the construction
from \autoref{eq:intersectfunctiongeneral} and \autoref{eq:functioningeneral}.

\begin{lemma} \label{lem:hyperplanesection}
  Let $H$ denote the standard hyperplane in $\R^N$. 
	For $i = 0, \dots, n$, we have 
	\[
	  H^{n-i} \cdot \Sigma_M = \Sigma_{M^{\leq i}}.
	\]
\end{lemma}

\begin{proof}
  By induction, it is sufficient to prove $H \cdot \Sigma_M = \Sigma_{M^{\leq n-1}}$.
	Setting 
	\[
	  h' = \max\{x_0, \dots, x_N\},
	\]
	we have $H \cdot \Sigma_M = (h' \cdot \Sigma'_M)/\R \mathbf{1}$,
	hence it remains to prove $h' \cdot \Sigma'_M = \Sigma'_{M^{\leq n-1}}$.
	By \autoref{lem:genericmatroids} and \autoref{eq:intersectfunctiongeneral}, 
	we can write $\Sigma'_{M^{\leq n-1}}$
	as $g'_1 \cdot \Sigma'_M$. By \autoref{eq:functioningeneral}, 
	$g'_1(S) = -1$ if $\rk(S) = n+1$ and $g'_1(S) = 0$ otherwise. For flats $S=F$ of $M$, $\rk(F) = n+1$ is equivalent to $F = E$. 
	On the other hand, we have $h'(S) = -1$ if and only if $S = E$. It follows that the functions $h'$ and $g'_1$ agree on
	flats of $M$ and hence on $\Sigma'_M$, which proves the claim. 
\end{proof}

\subsection{Cutting out matroid fans}

As a side remark, let us have a quick look at the opposite situation to the previous 
subsection, the inclusion $\Sigma_M \subset \R^N = \Sigma_{U_{N+1,N+1}}$.
In this case, \autoref{eq:intersectfunctiongeneral} allows us to write $\Sigma_M$ 
as a complete intersection
\begin{equation} \label{eq:matroidequations} 
  \Sigma_M = g_{N-n} \cdots g_1 \cdot \R^N.
\end{equation}
Of course, this is a complete intersection in a rather weak sense. For example,
the functions need not be tropically linear nor convex (tropical polynomials) in general.
Nevertheless, this description of arbitrary matroid fans might be useful in some contexts.
For example, the CSM classes of matroid fans studied in \cite{LRS-ChernSchwartzMacpherson}
can be written in terms of these functions as follows. 
Here, given a fan cycle $X \subset \R^N$ and a piecewise linear function $g$ on $X$, 
we consider $1+g$ as the operator on the group $Z_*(X)$ of fan cycles in $X$ given by
\[
  (1+g) \colon Z_*(X) \to Z_*(X), \;\;\; Z \mapsto Z + g \cdot Z. 
\]
The inverse of this operator is given by 
\[
  \frac{1}{1+g} = 1 - g + g^2 - g^3 + \cdots \colon Z_*(X) \to Z_*(X).
\]
(More formally, $1+g$ and $1/(1+g)$ are piecewise polynomial functions or tropical Chow cohomology elements, 
but this is not needed here).

\begin{proposition} 
  Let $\Sigma_M$ be a matroid fan cut out by the functions $g_i$ as in \autoref{eq:matroidequations}. 
	The the CSM-classes of $\Sigma_M$ are equal to 
  \[
    \text{CSM}_*(\Sigma_M) = \prod_{i=1}^{N-n} \frac{g_i}{1+g_i} \cdot \R^N = \prod_{i=1}^{N-n} \frac{1}{1+g_i} \cdot \Sigma_M.
	\]
\end{proposition}

\begin{proof}
  In the classical case, the analogous formula is essentially a consequence of the definitions and the adjunction formula.
	Here the quick and dirty proof goes as follows:  
	We proceed by induction on the codimension.
	The induction start $n=N$ is trivial, so let us assume $n < N$.
	Consider the generic chain $M = M_0, M_1, \dots, M_{N-n} = U_{N+1,N+1}$ associated to $\Sigma_M \subset \R^N$. 
	We set $\Sigma_0 = \Sigma_M = \Sigma_{M_0}$ and $\Sigma_1 = \Sigma_{M_1}$. 
	By induction assumption, we have $\text{CSM}_*(\Sigma_1) = \prod_{i=1}^{N-n-1} \frac{g_i}{1+g_i} \cdot \R^N$.
	We set $g = g_{N-n}$, hence $g \cdot \Sigma_1 = \Sigma_0$. 
	It remains to show that 
	\[
	  \text{CSM}_*(\Sigma_0) = \frac{g}{1+g} \cdot \text{CSM}_*(\Sigma_1).
	\]
	Since $M_0$ is an elementary quotient of $M_1$, there exists a unique matroid $L$ on the base set $E \sqcup \{e\}$ 
	such that $L \backslash e = M_1$ and $L / e = M_0$. Here, $e$ is some additional element
	which is neither loop nor coloop of $L$. We set $\Sigma = \Sigma_L \subset \R^{N+1}$.
	Let $\pi : \R^{N+1} \to \R^N$ denote the projection forgetting the
	coordinate $x_e$. 
	Then $\pi \colon \Sigma \to \Sigma_1$ is an elementary tropical modification along the function $g \colon \Sigma_1 \to \R$.
	By \cite[Proposition 5.4]{LRS-ChernSchwartzMacpherson}, we have
	\begin{equation} \label{eq:CSMModification} 
	  \text{CSM}_*(\Sigma) = \pi^* \text{CSM}_*(\Sigma_1) - \pi^* \text{CSM}_*(\Sigma_0).
	\end{equation}
	Given a fan cycle $X \in \R^{N+1}$, we denote by $D_e \cdot X$ the star fan $\Star_{-v_e}(X) \subset \R^N$ along the direction $-v_e$ 
	(we think of this geometrically as the intersection of $X$ with the divisor $D_e = \{x_e = -\infty\}$ at infinity).
	It follows from \cite[Lemma 8.8]{FR-DiagonalTropicalMatroid} and \cite[Proposition 5.3.13]{MR-TropicalGeometry} that 
	$D_e \cdot \pi^* X = g \cdot X$ for any fan cycle $X \subset \Sigma_1$. 
	Moreover, by the local definition of CSM-classes, it is clear that
	$D_e \cdot \text{CSM}_*(\Sigma) = \text{CSM}_*(D_e \cdot \Sigma) = \text{CSM}_*(\Sigma_0)$. 
	Hence, intersecting \autoref{eq:CSMModification} with $D_e$ we obtain 
	\[
	  \text{CSM}_*(\Sigma_0) = g \cdot (\text{CSM}_*(\Sigma_1) - \text{CSM}_*(\Sigma_2)). 
	\]
	Solving for $\text{CSM}_*(\Sigma_2)$, the result follows. 
\end{proof}

The $n_*$-classes appearing in \cite[Section 9.3]{LRS-ChernSchwartzMacpherson}
in connection with Speyer's $g$-polynomial \cite{Spe-MatroidInvariantVia}
can be written as 
\[
  n_*(\Sigma_M) = \frac{\prod_{i=1}^{N-n} (1+g_i)}{1+h} \cdot \Sigma_M.
\]
Here, $h = \max\{x_0, \dots, x_N\} - x_0$.
Finally, note that $g_1 = h$ if and only if $M$ has no coloops. 
Hence, in this case $n_*$ can be further simplified by clearing the
denominator. In view of the $g$-polynomial conjecture, it is
interesting to study the positivity properties
of the functions $g_i$ and the expressions above.

\section{Matroidal automorphisms and their beta invariant}

\subsection{Generalized beta invariants} \label{betainvariant}

Beta invariants are usually defined in the context of geometric lattices. 
We want to extend the definition to the case of arbitrary lattices $L$
equipped with a rank function $\rk$. 
The application we have in mind is the sublattice $L = \Fix(L(M))$ of fixed flats
of an matroidal automorphism (with
the restricted rank function), see \autoref{matautom}.
We define the generalized beta invariant as (up to signs) 
the convolution of $\rk$ with the Möbius function $\mu$ of $L$. 
The main property that we will need later (\autoref{lem:RecursiveBeta})
is based purely on this definition, without further requirements on $(L, \rk)$.

Let $L$ by a finite lattice with partial order $\subset$, 
minimal element $\emptyset$ and maximal element $E$. 
Let $\rk \colon L \to \Z$ be an arbitrary function on $L$, called the rank function.
We set $n := \rk(E) - \rk(\emptyset) -1$. 

\noindent
\emph{Convention}.
  We emphasize that in this text the symbol $\subset$ is used synonymously to $\subseteq$
	or $\leq$, that is, equality is included.

\begin{definition} 
  The \emph{beta invariant} of $L$ (or rather, $(L,\rk)$) is 
	\[
	  \beta(L) := (-1)^{n+1} \sum_{F \in L} \mu(\emptyset, F) \rk(F).
	\]
\end{definition}

For any $F \in L$, we equip the interval $[F,E]$ 
with the restricted rank function $\rk|_{[F,E]}$ and set
\begin{equation} \label{eq:betaforF} 
  \beta(F) := \beta([F,E]) = (-1)^{\rk(E) - \rk(F)} \sum_{F \subset G \in L} \mu(F, G) \rk(G).
\end{equation}

\begin{lemma} \label{lem:InverseBeta}
  For any $G \in L$, we have
	\[
	  \rk(G) = \sum_{G \subset F \in L} (-1)^{\rk(E)-\rk(F)} \beta(F).
	\]
\end{lemma}

\begin{proof}
  This is just Möbius inversion for \autoref{eq:betaforF}.
\end{proof}

For our purposes, it will be useful to rewrite this formula in a particular, asymmetric way.

\begin{lemma} \label{lem:RecursiveBeta}
  Fix an element $G \in L$, $G \neq \emptyset$. Then
	\[
	  (-1)^{n} \beta(L) = \rk(G) - \rk(\emptyset) \; - \; \sum_{\emptyset \neq F \notin [G,E]} (-1)^{\rk(E) -\rk(F)-1} \beta(F).
	\]
\end{lemma}

\begin{proof}
  Since $\beta(\emptyset) = \beta(L)$, 
	the formula is a rearrangement (including some sign yoga) of 
	\[
	  \rk(\emptyset) - \rk(G) = \sum_{F \notin [G,E]} (-1)^{\rk(E) -\rk(F)} \beta(F).
	\]
	This follows from using \autoref{lem:InverseBeta} twice.
\end{proof}

\begin{remark} 
  Let $K \subset L$ be a sublattice with $\emptyset, E \in K$. We equip 
  $K$ with the restricted rank function $\rk|_K$. 
	For any $S \subset L$, we set 
	\[
	  \cl(F) = \bigwedge_{\substack{G  \in K \\ F \subset G}} G.
	\]	
	By general properties of Möbius functions, the 
	Möbius function $\mu^K$ of $K$ can computed in terms
	of the Möbius function $\mu^L$ of $L$ by
	\[
	  \mu^K(\emptyset, G) = \sum_{\substack{F  \in L \\ \cl(F) = G}} \mu^L(\emptyset, F).
	\]
	It follows that the beta invariant of $K$ can be expressed as 
	\begin{equation} \label{eq:reformulationmob} 
		\beta(K) = (-1)^{n+1} \sum_{G \in K} \mu^K(\emptyset, G) \rk(G) = (-1)^{n+1} \sum_{F \in L} \mu^L(\emptyset, F) \rk(\cl(F)).
	\end{equation}
\end{remark}

\subsection{Matroidal automorphisms} \label{matautom}

A matroid automorphisms $\psi \colon M \to M$ is a bijection
$\psi \colon E \to E$ such that $F$ is a flat if and only if $\psi(F)$ is a flat. 
This induces an automorphism of geometric lattices $\psi \colon L(M) \to L(M)$
which we denote by the same letter. 
Vice versa, a lattice automorphism $\psi \colon L(M) \to L(M)$
defines a matroid automorphism 
up to the ambiguity of permuting parallel elements. 
For what follows, it is actually enough to fix the
lattice automorphism $\psi \colon L(M) \to L(M)$
(since the linear span of $\Sigma_M$ is generated
by the indicator vectors of flats).

\begin{definition} 
  Let $M$ be a loopless matroid. 
	Given an matroid automorphism $\psi \colon M \to M$, the associated
	\emph{matroidal automorphism} $\Psi \colon \Sigma_M \to \Sigma_M$ is the 
	restriction of the linear map on $\R^N$
	given by
	\[
	  \Psi \colon v_S \mapsto v_{\psi(S)}
	\]
	for all $S \subset E$. 
\end{definition}

It follows from the definition that $\Psi$ is a tropical automorphism of $\Sigma_M$ 
which respects its fine subdivision.

In the context of the trace formula, we are interested in the subset 
\[
  \Fix(L(M)) = \{F \in L(M) : F = \psi(F)\}.
\]

\begin{lemma} 
  The set $\Fix(L(M))$  is a sublattice of $L(M)$. 
\end{lemma}

\begin{proof}
  We need to show that $\Fix(L(M))$ is closed under $\wedge$ (intersection) and $\vee$ (union and closure).
	Since $\psi$ is a bijection, taking images commutes with intersection/unions, and since it takes
	flats to flats, it also commutes with taking closures. The statement follows. 
\end{proof}

\begin{remark} 
  As a sublattice of $L(M)$, $\Fix(L(M))$ is distributive. Note, however, 
	that in general it is neither graded, atomic, coatomic nor complemented. 
\end{remark}

In the context of \autoref{betainvariant}, 
we will always equip $\Fix(L(M))$ with the rank function $\rk$ of $L(M)$
restricted to $\Fix(L(M))$.
In other words, we set
\[
  \beta(\Fix(L(M))) = (-1)^{n+1} \sum_{\substack{F \in L(M) \\ F = \psi(F)}} \mu^\psi(\emptyset, F) \rk(F)
\]
and for any $F \in \Fix(L(M))$
\begin{equation} \label{eqbetainterval}
  \beta(\Fix([F,E])) = (-1)^{\rk(E) - \rk(F)} \sum_{\substack{G \in L(M) \\ F \subset G = \psi(G)}} \mu^\psi(F, G) \rk(G).
\end{equation}
Here, $\mu^\psi$ denotes the Möbius function on $\Fix(L(M))$ (in contrast to the Möbius function 
of $L(M)$). So, to emphasize, the definition mixes the 
Möbius function of $\Fix(L(M))$ with the rank function of $L(M)$. 
Note that given $F \in \Fix(L(M))$, we could alternatively consider the 
contracted matroid $M/F$ with rank function $\rk' = \rk - \rk(F)$ and induced automorphism $\psi'$. 
The beta invariant
\begin{equation}
  \beta(\Fix(M/F)) = (-1)^{\rk(E) - \rk(F)} \sum_{\substack{G \in L(M) \\ F \subset G = \psi(G)}} \mu^\psi(F, G) (\rk(G) - \rk(F)).
\end{equation}
is equal to $\beta(\Fix([F,E]))$ from \autoref{eqbetainterval}, except for $F=E$
(since summing the Möbius function over a non-trivial interval gives zero). 
In the (non-relevant) case $F=E$, we have $\beta(\{E\}) = \rk(E) = n+1$. 
Finally, note that 
\autoref{eq:reformulationmob} applied to $K = \Fix(L(M)) \subset L(M) = L$ 
allows us to rewrite the beta invariant as
\[
  \beta(\Fix(L(M))) = (-1)^{n+1} \sum_{F \in L(M)} \mu(\emptyset, F) \rk(F^\psi)
\]
where $F^\psi$ denotes smallest flat containing $F$ and fixed under $\psi$, 
\[
  F^\psi := \bigcap_{\substack{G \supset F \\ G = \psi(G)}} G.
\]

We recall from the introduction that our main theorem consists in the following
equations.
\[ 
	\deg (\Gamma_{\Psi} \cdot \Delta) = (-1)^n \beta(\Fix(L(M))) = 
		\sum_{p} (-1)^{p} \Tr(\Psi_*, \F_p(\Sigma_M))
\]
The next two sections are devoted to the proofs of the first and second equality, respectively.

\section{The intersection-theoretic side}

In this section, our goal is to prove
\begin{equation} \label{eq:rightside} 
  \deg(\Gamma_\Psi \cdot \Delta) = (-1)^n \beta(\Fix(L(M))).
\end{equation}

\subsection{General approach}

We denote by
$\Gamma \colon \Sigma_M \to \Sigma_M \times \Sigma_M$, $x \mapsto (x,\Psi(x))$
the graph map. 
We denote by $f_i := \Gamma^*(g_i)$ the pullbacks of the functions $g_i$
constructed in \autoref{eq:functiongidehom}. 

\begin{lemma} \label{lem:reformIntersection}
  With the notations from above, we have
	\[
	  \deg (\Gamma_\psi \cdot \Delta) = \deg(f_n \cdots f_1 \cdot \Sigma_M). 
	\]
\end{lemma}

\begin{proof}
  This follows from \autoref{completeintersection}, \cite[Theorem 4.5 (6)]{FR-DiagonalTropicalMatroid} and the 
	projection formula \cite[Proposition 7.7]{AR-FirstStepsTropical}.
\end{proof}

We set $X_k := f_k \cdots f_1 \cdot X$. Our goal is to compute 
information about these intermediate intersection products inductively.
The first new technical difficulty that we encounter in comparison
to \cite{Rau-TropicalPoincareHopf} is the fact 
that $\Gamma \colon \Sigma_M \to \Sigma_{M\otimes_0 M}$
is in general \emph{not} a map of fans, i.e.\ the cones of $\Sigma_M$ 
are not necessarily mapped to cones of $\Sigma_{M\oplus_0 M}$. 
Consequently, the functions $f_k$ are in general not linear when restricted
to cones $\sigma_\FF \subset \Sigma_M$. 
In principal, this could lead to fan structures which
cannot be described as subfans of the permutahedral fan.
However, somewhat mysteriously, it turns out that no such refinements are necessary. 
In fact, we will show that $f_k$ is linear when restricted to a (non-zero) cone
of $X_{k-1}$ (but not of $\Sigma_M$). By induction, this is enough to ensure that
$X_k$ can (still) be written as a subfan of the permutahedral fan. 
The following lemma contains the technical key piece of this argument.

\begin{lemma} \label{lem:fklinearity}
  Pick $k \in \{1, \dots, n\}$ and let $\FF = (F_i)$ be a chain of flats such that one of the following conditions holds:
	\begin{enumerate}
		\item For any $i$, either $0 \in F_i \cap \psi(F_i)$ or $0 \notin F_i \cup \psi(F_i)$. 
		\item The gap sequence of $\FF$ is $\gap(\FF) = (k-1, 0, \dots, 0)$. 
	\end{enumerate}
	Then the function $f_k$ is linear on $\sigma_\FF$. 
\end{lemma}

\begin{proof}
  If property (a) holds, then $\Gamma(\sigma_\FF)$ is a cone in the fine subdivision of
	$\Sigma_{M\otimes_0 M}$, given by the chain of flats $((F_i,\psi(F_i)))$.
	Since the $g_k$ are linear on such a cone by definition, the $f_k$ are linear
	on $\sigma_\FF$ as well.
	
	Let us now assume property (b) holds. We set $G_i = \psi(F_i)$ and 
	choose $p$ and $q$ maximal with the property that $0 \in F_p$ and $0 \in G_q$, respectively. 
	If $p=q$, we are back in the case of property (a). By symmetry of the function $g_k$ 
	in the two factors, we can assume $p > q$ without loss of generality. 
	Given a point $x \in \sigma_\FF$, we write it as a positive linear combination
	of the $v_{F_i}$ with coefficients $a_i$. 
	We claim that $f_k|_{\sigma_\FF}(x)$ 
	is equal to the linear function $A(x) := \sum_{i=1}^p a_i$. 
	
	Set $y := \Gamma(x) = (x, \psi(x))$. Using indicator vectors, we can write $y$
	as 
	\begin{equation} \begin{split} \label{eq:vectorsum} 
	  y &= a_1 v_{(F_1, G_1)} + \dots + a_q v_{(F_q, G_q)} \\
		  &+ a_{q+1} (v_{(F_{q+1}, E)} + v_{(\emptyset, G_{q+1})}) + \dots + a_p (v_{(F_{p}, E)} + v_{(\emptyset, G_{p})}) \\
			&+	a_{p+1} v_{(F_{p+1}, G_{p+1})} + \dots + a_l v_{(F_l, G_l)}.
	\end{split} \end{equation}
	Here, we use the formula $(v_F, v_G) = v_{(F, E)} + v_{(\emptyset, G)}$ if $0 \in F\setminus G$,
	which follows from our conventions.
	Note that when considered as an equation in $\R^{E'} = \R^{2N+1}$, $y$ is normalised
	in the sense that its minimal coordinate is $0$. Moreover, $y_0 = A(x)$. 
	More generally, the coefficient $y_e$ of $y$,
	depending on whether $e$ is chosen from the first or second copy of $E$, 
	is equal to $b_i := a_1 + \dots + a_i$
	or $c_i := (a_{q+1} + \dots a_p) + b_i$ where $i \in \{0, \dots, l\}$
	denotes the maximal $i$ such that $e \in F_i$ or $e \in G_i$,
	respectively. 
	It follows that the cone 
	$\sigma \subset \Sigma_{M \otimes_0 M}$ that contains $y$ in its relative interior
	corresponds to a chain of flats of the form $(F_i, G_j)$.
	Indeed, given $t \in \R$ denote by $i(t)$ and $j(t)$ the maximal 
	indices such that $b_{i(t)} \leq t$ and $c_{j(t)} \leq t$. 
	As we increase $t$, the tuples $(F_{i(t)}, G_{j(t)})$
	yield the required chain. Note that $b_p = c_q$, so indeed
	each of the tuples $(F_i, G_j)$ in this chain 
	satisfies $0 \in F_i \cap G_j$ or $0 \notin F_i \cup G_j$, 
	which can also be rewritten as either $i \leq p, j \leq q$ or $i > p, j > q$. 
	Note that $p \geq 1$, so in the latter case $i \geq 2$. 
	
	By our assumption, we have $\rk(F_i)=\rk(G_i) \leq n+1-k$ for $i \geq 1$ and
	$\rk(F_i)=\rk(G_i) < n-k+1$ for $i \geq 2$. 
	It follows that $\rk(F_i) + \rk(G_j) \leq \rk(F_i \cup G_j) + n+1-k$
	if $(i,j) \neq (0,0)$ and the inequality is strict 
	if $i > p, j > q$.
	Using the description of $g_k$ in \autoref{eq:functiongidehom}, 
	we find
	\[
	  g_k(F_i, G_j) = \begin{cases} 
		                  1 & i \leq p, j \leq q, (i,j) \neq (0,0), \\
		                  0 & i > p, j > q \text{ or } (i,j) = (0,0).
										\end{cases}
	\]
	We now write $y$ as a positive sum of the cone generators of $\sigma$. 
	Setting the coefficient $v_{(E,E)}$ to zero, this sum
	is normalised as above (minimal coordinate is $0$) and hence
	equal to \autoref{eq:vectorsum} in $\R^{2N+1}$.
	Moreover, by the previous computation, $g_k(y)$ (and hence $f_k(x)$)
	is equal to the sum of coefficients of the vectors
	$v_{(F_i, G_j)}$ with $0 \in F_i \cap G_j$, which clearly is equal to
	$y_0$. Since $y_0 = A(x)$, this proves the claim.
\end{proof}

\subsection{Combinatorial properties of $X_k$}

We will now proceed by formulating certain properties of the intermediate intersection
products $X_k$. 
We split the properties into two statements: Here we consider combinatorial properties, in the next subsection
we discuss weights.

\begin{lemma} \label{Xkcombinatorics}
  For all $k \in \{0, \dots, n\}$, the following statements hold.
	\begin{enumerate}
		\item \label{lab1}
		  The tropical cycle $X_k$ can be represented as a weighted subfan of the 
			permutahedral fan. (By a \emph{facet} of $X_k$, we mean a cone of non-zero weight in $X_k$.)
		\item \label{lab2}
		  For a facet $\sigma_\FF$ of $X_k$, the gap sequence of $\FF$ has one of the following two forms.
			\begin{align} 
				\gap(\FF) &= (r,s,0, \dots, 0) =: (A), & r+s=k, \\
				\gap(\FF) &= (r,s,0, \dots, 0,1,0, \dots, 0) =: (B), & r+s = k-1.
			\end{align}
		\item \label{lab3}
		  If $\gap(\FF) = (A)$, $0 \notin F_1 \cup \psi(F_1)$ and $s \geq 1$, then $F_1 = \psi(F_1)$.
		\item \label{lab4}
		  If $\gap(\FF) = (A)$, $0 \in F_1 \cup \psi(F_1)$ and $s \geq 1$, then $0$ and $\psi^{-1}(0)$ are linearly independent in $F_1/F_2$.
			If moreover $s \geq 2$, then $F_1 = \psi(F_1)$.
		\item \label{lab5}
		  If $\gap(\FF) = (B)$ and $G \supsetneq H$ denotes the part of $\FF$ corresponding
			to the final $1$ in $\gap(\FF)$,
			then $0$ and $\psi^{-1}(0)$ are linearly independent in $G/H$
			(and hence $G = \cl(H \cup\{ 0 , \psi^{-1}(0)\})$). 
			If moreover $s \geq 1$, then $F_1 = \psi(F_1)$.
		\item \label{lab6}
		  For $k < n$, $f_{k+1}$ is linear when restricted
			to a facet of $X_k$. 
	\end{enumerate}
\end{lemma}

The proof of this lemma mostly consists of simple, but tedious calculations
of weights in the intersection product $f_{k+1} \cdot X_k$. 
To make it as transparent as possible, we will first
collect a couple of recurrent arguments and case distinctions
that occur in this calculation.

\begin{remark} \label{rem:weightcomp}
  Let $\tau$ be a codimension one cone in $X_k$ and let $\FF$ be the corresponding 
	chain of flats. Our general goal is to compute the weight $\omega(\tau)$ of $\tau$
	in $f_{k+1} \cdot X_k$. 
  \begin{enumerate}
		\item 
		  The facets in $X_k$ containing $\tau$ correspond to filling a (non-trivial) gap 
			$G \supsetneq H$ of $\FF$ with an additional flat $G \supsetneq F \supsetneq H$.
			The balancing condition around $\tau$ naturally splits into separate equations,
			one for each gap of $\FF$ (and the facets corresponding to it). 
			In particular, the calculation of $\omega(\tau)$ can be split into a calculation
			for each gap.
		\item 
		  Let $F^1, \dots, F^m$ denote the flats corresponding to the facets 
			of $X_k$ for a given gap $G \supsetneq H$ of $\FF$. The typical situation
			will be that the sets $F^i \setminus H$ form a partition of $G$. In this situation,
			the involved indicator vectors satisfy a unique linear relation (up to multiples), namely
			\[
			  \sum_{i=1}^m v_{F^i} = v_G + (m-1) v_H. 
			\]
			By uniqueness, it follows that the weights of the corresponding facets in $X_k$ 
			are all equal, say, to $\omega \in \Z\setminus\{0\}$. 
		\item
		  By induction, we may assume that $f_{k+1}$ is linear on the facets
			of $X_k$ (which are cones of the permutahedral fan). It follows that
			in order to compute $\omega(\tau)$, it is sufficient to know
			the values $f_{k+1}(F)$. 
			Using \autoref{eq:functiongidehom}, these values are given by
			\begin{equation}  
			  f_{k+1}(F) = \begin{cases}
	                 -1 & 0 \notin F \cup \psi(F) \text{ and } 2 \rk(F) \geq \rk(F \cup \psi(F)) + n - k, \\
							     +1 & 0 \in F \cap \psi(F) \text{ and } 2 \rk(F) \leq \rk(F \cup \psi(F)) + n - k, \\
							     +1 & 0 \in (F \cup \psi(F)) \setminus (F \cap \psi(F)) \text{ and } \rk(F) \leq n - k, \\
						        0 & \text{otherwise}.
					       \end{cases}
			\end{equation}
			In the third case, we use the fact that $\Gamma(v_F) = v_{(F,E)} + v_{(\emptyset, \psi(F))}$ or
			$\Gamma(v_F) = v_{(F,\emptyset)} + v_{(E, \psi(F))}$, depending on whether $0 \in F$ or
			$0 \in \psi(F)$. It is convenient to list a few particular cases, focusing on the critical rank $n-k$.
			\begin{equation} \label{eq:functionsfi} 
			  f_{k+1}(F) = \begin{cases}
	                 -1 & 0 \notin F = \psi(F) \text{ and } \rk(F) \geq n - k, \\
	                  0 & 0 \notin F \cup \psi(F), F \neq \psi(F) \text{ and } \rk(F) = n - k, \\
	                  0 & 0 \notin F \cup \psi(F) \text{ and } \rk(F) < n - k, \\
										0 & 0 \in F = \psi(F) \text{ and } \rk(F) \geq n - k + 1, \\
									 +1 & 0 \in F \cap \psi(F), F \neq \psi(F) \text{ and } \rk(F) = n - k + 1, \\
							     +1 & 0 \in F \cap \psi(F) \text{ and } \rk(F) \leq n - k. \\
					       \end{cases}
			\end{equation}
		\item
		  Going back to the partition case of item (b), let $q$ be the number of flats
			$F^i$ such that $0 \in F^i \cup \psi(F^i)$. The possible values are
			$q = 0, 1, 2, m$. The first and latter case correspond to 
			$0 \notin G \cup \psi(G)$ and $0 \in H \cup \psi(H)$, respectively. 
			The case $q=2$ occurs if $0$ and $\psi^{-1}(0)$ are linearly independent
			in $G/H$. The remaining cases correspond to $q=1$. 
			Based on all the previous comments, we list in \autoref{weightcomputation}
			the computation of $\omega(\tau)$ (or rather, the contribution of
			a fixed gap $G \supsetneq H$ to it) for the various values of $q$
			and with various extra conditions.
	\end{enumerate}
\end{remark}

\begin{table}[tb]%
\centering
\begin{tabular}{l|l|l|l|ll}
\hline
               & $q=0$                          & $q=1$            & $q=2$            & $q=m$                          \\ \hline
condition      & $\rk G < n-k$                  & $\rk G \leq n-k$ & $\rk G \leq n-k$ & $\rk G \leq n-k$               \\
$f_{k+1}(G)$   & $0$                            & $1$              & $1$              & $1$                            \\
$f_{k+1}(F^i)$ & $0, \dots, 0$                  & $1, 0, \dots, 0$ & $1,1,0 \dots,0$  & $1, \dots, 1$                  \\
$f_{k+1}(H)$   & $0$                            & $0$              & $0$              & $1$                            \\
$\omega(\tau)$ & $0$                            & $0$              & $\omega$         & $0$                            \\ \hline
condition      & $\rk G = n-k$                  & \begin{tabular}[t]{@{}l@{}} 
                                                    $\rk G = n-k+1$, \\ 
																										$[0 \notin F^i \cup \psi(F^i)$ \\ 
																										$\Rightarrow \rk F^i < n-k]$ 
																									\end{tabular}    &                  & $\rk G = n-k+1$                \\
subcase        & $G =$/$\neq \psi(G)$           & $G =$/$\neq \psi(G)$ &              & $G =$/$\neq \psi(G)$           \\
$f_{k+1}(G)$   & $-1$/$0$                       & $0$/$1$          &                  & $0$/$1$                        \\
$f_{k+1}(F^i)$ & $0, \dots, 0$                  & $1, 0, \dots, 0$ &                  & $1, \dots, 1$                  \\
$f_{k+1}(H)$   & $0$                            & $0$              &                  & $1$                            \\
$\omega(\tau)$ & $\omega$/$0$                   & $\omega$/$0$     &                  & $\omega$/$0$                   \\ \hline
condition      & $G = \psi(G)$,                 &                  &                  & $G = E$,                       \\
               & $\rk F^i = n-k-1$              &                  &                  & $\rk F^i = n-k$                \\
$f_{k+1}(G)$   & $-1$                           &                  &                  & $0$                            \\
$f_{k+1}(F^i)$ & $0, \dots, 0$                  &                  &                  & $1,\dots,1$                    \\
$f_{k+1}(H)$   & $0$                            &                  &                  & $1$                            \\
$\omega(\tau)$ & $\omega$                       &                  &                  & $\omega$                       \\ \hline
\end{tabular}
\caption{The computation of the weight $\omega(\tau)$ for various types of facets $\tau$ in $f_{k+1} \cdot X_k$. 
The notation is borrowed from \autoref{rem:weightcomp}.}
\label{weightcomputation}
\end{table}

\begin{proof}[of \autoref{Xkcombinatorics}]
  The initial case $k = 0$ is trivial (except for 
	\autoref{lab6}, which follows by the same argument
	as below).
	
	Let us consider the induction step $k \to k+1$. 
	For \autoref{lab1}, note that by induction assumption 
	$X_{k}$ can be represented on the permutahedral fan
	and $f_{k+1}$ is linear on the facets of this representation.
	Hence $X_{k+1} = f_{k+1} \cdot X_{k}$ can also be represented on the
	permutahedral fan. 
	
	For \autoref{lab2}, let $\tau$ be a cone of dimension $n-k-1$ whose
	weight $\omega(\tau)$ in $X_{k+1}$ is non-zero. 
	Let $\FF$ be the corresponding chain of flats with gap sequence 
	$\gap(\FF) = (r_0, \dots, r_l)$. We need to show that
	$S := \sum_{i=2}^l r_i \leq 1 $. 
	If $S > 2$, then by induction assumption $\tau$ is not contained
	in any facet of $X_k$, a contradiction. 
	
	So let us assume $S=2$. 
	Note that hence $r_1 + r_2 = k-1$ and therefore
	$\rk F_2 = n-k$. If $(r_2, \dots, r_l) = (\dots, 1, \dots, 1, \dots)$
	(the dots represent zeros), one of the two gaps must be as 
	in \autoref{lab5} and the facets of $X_k$ correspond to filling the other
	gap. Depending on the ordering of the gaps, $\omega(\tau)$
	is computed according to \autoref{weightcomputation}, row 1, $q=0$ or $q=m$.
	In both cases $\omega(\tau) = 0$. 
	If $(r_2, \dots, r_l) = (\dots, 2, \dots)$, with gap $G \supsetneq H$, 
	then $0$ and $\psi^{-1}(0)$ must be linearly independent in $G/H$ 
	and the possible fillings are given by $F = \cl(H \cup\{ 0 , \psi^{-1}(0)\})$
	and $F \nsubset \cl(H \cup\{ 0 , \psi^{-1}(0)\})$, $\rk F = \rk H +1$. 
	Note that this still induces a partition of $G \setminus H$ 
	and $q=1$, so by \autoref{weightcomputation}, row 1, $q=1$ we get
	$\omega(\tau) = 0$ again. This finishes \autoref{lab2}.
	
	We proceed with \autoref{lab3}, so $\gap(\FF) = (r,s, \dots)$ with $r+s = k+1$,
	$s \geq 1$ and $0 \notin F_1 \cup \psi(F_1)$. We want to show
	$F_1 = \psi(F_1)$.
	Note that the only possible (non-zero) fillings must have gap sequence
	$(r,s-1, \dots)$ (since type $(B)$ is excluded by $0 \notin F_1 \cup \psi(F_1)$ 
	and \autoref{lab5}).
	For $s > 1$, the statement follows by induction assumption. For $s=1$,
	note that $r = k$ and hence $\rk F_1 = n - k$. So the statement follows 
	from \autoref{weightcomputation}, row 2, $q=0$. 
	
	For \autoref{lab4}, we assume $\gap(\FF) = (r,s, \dots)$ with $r+s = k+1$,
	$s \geq 1$ and $0 \in F_1 \cup \psi(F_1)$. Assume that $0$ and $\psi^{-1}(0)$
	are not linearly independent in $F_1/F_2$. Then we must have $s=1$ 
	and the only possible gap sequence for fillings is $(r,0, \dots)$
	(all other possibilities have no facets adjacent). 
	In this case $\rk F_1 = n-k$ and \autoref{weightcomputation}, row 1, 
	$q=1$ or $q=m$ gives $\omega(\tau) = 0$, a contradiction. 
	
	Now assume $s \geq 2$. We need to show $F_1 = \psi(F_1)$. 
	The only possible gap sequences of fillings are 
	$(r,s-1,0,\dots)$ and $(r,s-2,1, \dots)$. 
	For $s > 2$, the statement follows
	by the induction assumption. For $s=2$, we have $\rk(F_1) = n-k+1$ and
	the possible fillings agree with the ones described in the last case of \autoref{lab1}.
	This is covered by \autoref{weightcomputation}, row 2, $q=1$.
	
	Let us now consider \autoref{lab5}, so $\gap(\FF) = (r,s, \dots, 1, \dots)$.
	Assume first that $0$ and $\psi^{-1}(0)$ are not linearly independent in $G \supsetneq H$
	(the final gap in $\FF$). Then the only possible fillings are fillings of $G \supsetneq H$. 
	But the corresponding weight can be computed according to \autoref{weightcomputation}, 
	row 1, and we get a non-zero weight only if $q=2$, a contradiction. 
	
	Now let $s \geq 1$ and assume that $F_1 \neq \psi(F_1)$. 
	Since by the previous argument $0 \in F_2 \cap \psi(F_2)$, the only possible 
	fillings have gap sequence $(r,s-1, \dots, 1, \dots)$. Again, if $s > 1$, 
	the statement follows from the induction assumption. If $s=1$, we 
	have $\rk F_1 = n-k+1$ and hence the claim follows from \autoref{weightcomputation}, 
	row 2, $q=m$.
	This finishes the proof of \autoref{lab5}.
	
	Finally, let us prove \autoref{lab6}. Note that up to now, we established that
	$X_{k+1}$ can be represented as a weighted subfan of the permutahedral fan and
	that its facets satisfy the properties of \autoref{lab3} -- \autoref{lab5}.
	Then the linearity of $f_{k+2}$ follows from \autoref{lem:fklinearity}. Indeed,
	note that all facets of $X_{k+1}$ satisfy condition (a) of that lemma, except for 
	the facets from \autoref{lab3} or \autoref{lab4} with $s=0$. 
	These ones satisfy condition (b) instead.
\end{proof}

\subsection{Weights on $X_k$}

Based on our understanding of the combinatorics of $X_k$, we can now describe
the weights of some of its facets.

\begin{lemma} \label{lem:weightsXk}
  Fix $k \in \{0, \dots, n\}$ and let $\sigma = \sigma_\FF$ be a facet of $X_k$. 
	Then the following statements hold.
	\begin{enumerate}
		\item \label{labb1}
		  If $\gap(\FF) = (k,0, \dots, 0)$ and $0 \in F_1 \cup \psi(F_1)$, then 
			$\omega(\sigma) = 1$.
		\item \label{labb2}
		  If $\gap(\FF) = (k-1,0, \dots, 0,1,0, \dots, 0)$ and $F_1 \neq \psi(F_1)$, then
			$\omega(\sigma) = \rk F_1^\psi - \rk F_1$.
		\item \label{labb3}
		  If $0 \notin F_1 $ and $F_1 = \psi(F_1)$, then 
			$\omega(\sigma) = (-1)^{n-\rk F_1} \beta(\Fix(M/F_1))$.
			
	\end{enumerate}
\end{lemma} 

Before proving the lemma, let us check that 
it implies \autoref{eq:rightside} as promised. 

\begin{proof}[\autoref{eq:rightside}]
  In the case $k=n$, the only chain of correct dimension is the 
	trivial flag $\FF = (E \supset \emptyset)$.
	We have $\gap(\FF) = (n)$ and $0 \notin F_1 \cup \psi(F_1) = \emptyset$.
	Hence, by item (c) of \autoref{lem:weightsXk} and \autoref{lem:reformIntersection},
	we conclude
	\[
	  \deg(\Gamma_\Psi \cdot \Delta) = \deg(X_n) = \omega(\sigma_{(E \supset \emptyset)}) = (-1)^n \beta(\Fix^\psi(M)).
	\]
\end{proof}

We now want to prove \autoref{lem:weightsXk}. 

\begin{proof}[\autoref{lem:weightsXk}]
  We (again) proceed by induction on $k$. For $k=0$, the statements are trivial (only \autoref{labb1} and
	\autoref{labb3} occur, and both give weight $1$ in this case).
	
	We now prove the induction step $k \to k+1$. Let $\sigma = \sigma_\FF$ be a facet of $X_{k+1}$. 
	We start with \autoref{labb1}. In this case, the facets of $X_k$ containing $\sigma$
	correspond to fillings $E \supsetneq F \supsetneq F_1$ with gap sequences of the form $(r,s, \dots)$,
	$r+s = k$. 
	Note that since $0 \in F_1 \cup \psi(F_1)$ (now the second step of the chain), 
	by \autoref{Xkcombinatorics} only the value $s=0$ is possible. Therefore, by induction
	assumption, all facets containing $\sigma$ have weight $1$, and the computation
	for $\omega(\sigma) = 1$ is given in \autoref{weightcomputation}, row 3, $q=m$. 
	
	We proceed with \autoref{labb2}. In this case, we have two gaps that can potentially be filled, namely
	$E \supsetneq F_1$ and $G \supsetneq H$, the gap corresponding to the final $1$. 
	By \autoref{Xkcombinatorics}, $0$ and $\psi^{-1}(0)$ are independent in $G/H$. 
	Hence the contribution $\omega_1$ of the fillings of $G \supsetneq H$ to $\omega(\sigma)$ can be computed according to 
	\autoref{weightcomputation}, row 1, $q=2$. Note that these fillings	correspond to facets of $X_k$ of 
	the type discussed in \autoref{labb1}, hence by induction assumption the all have weight $1$. 
	It follows that $\omega_1 = 1$. 
	
	Let us now consider the gap $E \supsetneq F_1$. 
	This is the first time that we encounter a facet structure
	that is not just given by a partition of $E \setminus F_1$. In fact, by \autoref{Xkcombinatorics}, the possible 
	fillings are given by flats $F \supsetneq F_1$ which satisfy at least one of the following conditions: 
	Either $\rk F = \rk F_1 + 1 = n-k+1$, or $F = \psi(F)$. Denoting the weights of the corresponding facets in $X_k$ by
	$\omega(F)$, the balancing condition for $X_k$ states that there are coefficients $\omega(E), \omega(F_1) \in \Z$
	such that 
	\begin{equation} \label{eq:balancingSpecial} 
	  \sum_F \omega(F) v_F = \omega(E) v_E + \omega(F_1) v_{F_1}.
	\end{equation}
	Pick an element $i \in E \setminus F_1$. Then we can express the coefficients $\omega(E), \omega(F_1)$
	as 
	\begin{equation} \begin{split} 
		\omega(E) &= \sum_{F \ni i} \omega(F), \\
		\omega(F_1) &= \sum_F \omega(F) \; - \omega(E) = \sum_{F \not\ni i} \omega(F). 
	\end{split} \end{equation}
	By assumption, $F_1 \neq F_1^\psi$ and we can choose $i \in F_1^\psi \setminus F_1$. 
	This implies $i \in F$ for all $F = \psi(F)$, so for this choice of $i$ 
	the sum for $\omega(F_1)$ can be restricted 
	to $F$ with  $F \neq \psi(F)$ (and hence $\rk F = \rk F_1 + 1$). 
	Finally, by comparing with \autoref{eq:functionsfi}, we find that the 
	values of the indicator vectors in 	\autoref{eq:balancingSpecial} 
	under $f_{k+1}$ are zero except for the cases 
	$F_1$ and $F \neq \psi(F)$, in which case the value is $1$.
	We conclude that the contribution $\omega_2$ of $E \supsetneq F_1$ is 
  \[
	  \omega_2 = \sum_{\substack{F \\ F \neq \psi(F)}} \omega(F) \; - \omega(F_1) 
		  = \sum_{\substack{F \\ F \neq \psi(F)}} \omega(F) - \sum_{\substack{F \not\ni i \\ F \neq \psi(F)}} \omega(F) 
			= \begin{cases} 
			    \omega(F') & \text{if } F' \neq F_1^\psi \\
					0 & \text{if } F' = F_1^\psi
				\end{cases}
	\]
	using the shorthand $F' = \cl(F_1 \cup \{i\})$.
	Note that we have $(F')^\psi = F_1^\psi$, so either by 
	induction assumption (if $F' \neq F_1^\psi$)
	or trivially (if $F' = F_1^\psi$)
	we have 
	\[
	  \omega_2 = \rk F_1^\psi - \rk \cl(F_1 \cup \{i\}) = \rk F_1^\psi - \rk F_1 - 1.
	\]
	So finally we get $\omega(\sigma) = \omega_1 + \omega_2 = \rk F_1^\psi - \rk F_1$, which proves \autoref{labb2}.
	
	It remains to prove \autoref{labb3}. In this case $0 \notin F_1 \cup \psi(F_1)$ and hence
	$\gap(\FF) = (r,s, \dots)$, $r+s = k+1$. Assume first that $s \geq 1$. Then by \autoref{Xkcombinatorics}
	the only possible fillings have gap sequence $(r, s-1, \dots)$. Note that $\rk F_2 = n-k-2$ while
	the ranks of the fillings are $\rk F = \rk F_2 + 1 = n-k-1$. We can thus use 
	\autoref{weightcomputation}, row 3, $q=0$, which proves the claim in this case.
	
	Assume now $s=0$, so $\gap(\FF) = (k+1, \dots)$. This is another case
	where the balancing condition cannot be expressed in terms of a partition of $E \setminus F_1$.
	Note that $\rk F_1 = n-k-1$ and hence
	$f_{k+1}(F_1) = f_{k+1}(E) = 0$. So, in order to compute $\omega(\sigma)$, we are only 
	interested in fillings $E \supsetneq F \supsetneq F_1$ which correspond to facets of $X_k$, on the
	one hand, and for which $f_{k+1}(F) \neq 0$, on the other hand. By comparing \autoref{eq:functionsfi}
	and \autoref{Xkcombinatorics}, we see that such flats $F$ belong to one of the following three subcases:
	\begin{enumerate}
		\item[(i)] $0 \notin F$, $F = \psi(F)$
		\item[(ii)] $0 \in F \cup \psi(F)$, $\rk F = n-k$
		\item[(iii)] $0 \in F \cup \psi(F)$, $\rk F = n-k+1$, $F = \cl(F_1 \cup \{0, \psi^{-1}(0)\})$, $F \neq \psi(F)$
	\end{enumerate}
	We set $G = \cl(F_1 \cup \{0\})$. 
	Note that $f_{k+1}(F) = -1$ for item (i) and $f_{k+1}(F) = 1$ 
	for items (ii) and (iii). 
	By the induction assumption, the flats of type (i) account for the last term in
	\autoref{lem:RecursiveBeta} applied to $L= \Fix^\psi(M/F_1)$ and $G^\psi$.
	Hence it remains to show that the contribution $\omega'$ of items (ii) and (iii) to $\omega(\sigma)$
	is $\rk(G^\psi)-\rk(F_1)$ (the first two terms in \autoref{lem:RecursiveBeta}).
	This can easily be checked by going through the following case consideration.
	\begin{itemize}
		\item $\rk(G^\psi) = \rk(F_1) + 1$ --- In this case, $F= G^\psi = G$ is the only flag of type (ii) and type (iii)
			does not occur, so $\omega' = 1$. 
		\item $\rk(G^\psi) = \rk(F_1) + 2$ --- In this case, $G$ and $\cl(F_1 \cup \{\psi^{-1}(0)\})$ are distinct and
      contribute to (ii), but $\cl(F_1 \cup \{0, \psi^{-1}(0)\}) \in \Fix^\psi(L(M))$, so again (iii) does not occur, so $\omega' = 2$.
		\item $\rk(G^\psi) > \rk(F_1) + 2$ --- Again, $G$ and $\cl(F_1 \cup \{\psi^{-1}(0)\})$ are distinct and
      contribute to (ii). Moreover, $F = \cl(F_1 \cup \{0, \psi^{-1}(0)\}) \notin \Fix^\psi(L(M))$, so it contibrutes
			to (iii). By induction assumption, the weight of the corresponding facet is 
			\[
			  \rk F^\psi - \rk F = \rk G^\psi - (\rk F_1 + 2) = \rk(G^\psi) - \rk(F_1) - 2.
			\]
			Hence $\omega' = 2 + (\rk(G^\psi) - \rk(F_1) - 2) = \rk(G^\psi) - \rk(F_1)$.
	\end{itemize}
	This finishes the proof of the lemma. 
\end{proof}

\begin{remark} 
  In contrast to \cite{Rau-TropicalPoincareHopf}, we are unable to give
	a complete description of the intermediate cycles $X_k$. 
	Rather, the properties established in 
	\autoref{Xkcombinatorics} and \autoref{lem:weightsXk}
	are chosen such as to be sufficient 
	to make the induction run, on the one hand,
	and to prove the main statement for $k=n$, on the other hand. 
	Nevertheless, properties such as \autoref{lem:weightsXk}
	\autoref{labb3} indicate that the $X_k$ exhibit
	an interesting recursive structure, similar 
	to the analogous cycles in \cite[Section 3]{Rau-TropicalPoincareHopf}
	The appearance of certain
	beta invariants is reminiscent of the so-called CSM classes 
	defined in \cite{LRS-ChernSchwartzMacpherson}.
	More general, a future goal is to find expressions of
	the cycles $X_k$ in terms of (pullbacks/pushforwards along $\Psi$ of) 
	the various 
	canonical classes on matroid fans defined
	in \cite{FS-kClassesMatroids,AHK-HodgeTheoryCombinatorial,
	LRS-ChernSchwartzMacpherson,ADH-LagrangianGeometryMatroids, 
  AP-HodgeTheoryTropical,BEST-TautologicalClassesMatroids}.
\end{remark}

\section{The trace side}

In this section we prove the second half of the main theorem by showing 
\begin{equation} \label{eq:secondhalf} 
	(-1)^n \beta(\Fix(L(M))) = 
		\sum_{p} (-1)^{p} \Tr(\Psi_*, \F_p(\Sigma_M)).
\end{equation}

\subsection{A resolution of the framing groups}

We start by constructing an explicit resolution for the 
framing groups $\F_p(\Sigma_M)$ of a matroid fan using chains of flats of low rank.
It appears that this resolution is known to experts --- a similar resolution is for example 
considered in \cite{AP-HodgeTheoryTropical}. We give a brief independent treatment here. 

Let $\FF = E \supsetneq F_1 \supsetneq \dots \supsetneq F_l \supsetneq  \emptyset$ 
be a chain of flats of $M$. 
Recall that $l(\FF) := l$ denotes the \emph{length} of $\FF$.
Moreover, we define the \emph{rank} of $\FF$ by $\rk(\FF) := \rk(F_1)$. 
Hence, except for $E$, $\FF$ only involves flats of rank at most $\rk(\FF)$. 
Let $\CC^{\leq p}_l$ be the set of chains $\FF$ of length $l$ and rank at most $p$ and let $\R  \CC^{\leq p}_l$
be the real vectorspace with a basis labelled by $\CC^{\leq p}_l$.

We consider the differential complex formed by the simplicial coboundary maps 
$\partial \colon \R \CC^{\leq p}_l \to \R \CC^{\leq p}_{l+1}$.
Explicitly, $\partial$ maps a generator $e_\FF, \FF \in \CC^{\leq p}_l$, to a vector
whose non-zero entries correspond to the chains $\GG \in \CC^{\leq p}_{l+1}$ with 
$\FF \leq \GG$. Moreover, such an entry is equal to $(-1)^k$ where $k$ denotes the index
such that $G_k \notin \FF$. 
Given a chain $\FF$ of length $p$, we set 
\[
  V_\FF := v_{F_1} \wedge \dots \wedge v_{F_p} \hspace{3ex} \in \F_p(\Sigma_M),
\]
or, in other words, $V_\FF$ is the canonical volume element for the cone $\sigma_\FF \subset \Sigma_M$.
This defines, by construction of $\F_p(\Sigma_M)$, a surjective map $\R \CC_p \to \F_p(\Sigma_M)$.
We are interested in its restriction $\R  \CC^{\leq p}_p \to \F_p(\Sigma_M)$ to chains 
only using flats of rank at most $p$. Our goal ist prove the following statement.

\begin{theorem} \label{thm:res}
  Given a matroid fan $\Sigma_M$ with framing groups $\F_p(\Sigma_M)$, the sequence 
	\[
	  0 \to \R \CC^{\leq p}_0 \to \dots \to \R \CC^{\leq p}_p \to  \F_p(\Sigma_M) \to 0
	\]
	is exact for all $p$.
\end{theorem}

We start by proving the theorem for $p = n = \dim(\Sigma_M)$. In a sceond step, we show how
to reduce to this case by using a tropical analogue of the Lefschetz hyperplane section theorem.

\begin{lemma} \label{lem:resmaxdim}
  Given a matroid fan $\Sigma_M$ of dimension $n$, the sequence 
	\begin{equation} \label{resmaxdim} 
	  0 \to \R \CC_0 \to \dots \to \R \CC_n \to  \F_n(\Sigma_M) \to 0
	\end{equation}
	is exact.
\end{lemma}

\begin{proof}
  Recall that by Poincaré duality \cite{JSS-SuperformsTropicalCohomology}
	we have 
	\[
	  \F_n(\Sigma_M) = H_0(\Sigma_M, \F_n) \cong (H_n^{\BM}(X, \Z))^*.
	\]
	Moreover, note that $H_k^{\BM}(X, \Z) = 0$ for all $k \neq 0$ (again, by Poincaré duality, or
	using the well-known statement about homology groups of geometric lattices \cite{Fol-HomologyGroupsLattice}).
	Hence the complex of simplicial chains, completed by $H_n^{\BM}(X, \Z)$,
	\begin{equation} 
		0 \to H_n^{\BM}(X, \Z) \to \R \CC_n \to \dots \to \R \CC_0 \to 0
	\end{equation}
	is exact and is dual to \autoref{resmaxdim} under this identification. 
\end{proof}

\begin{remark} 
  Since Poincaré duality holds over $\Z$ by \cite{JRS-Lefschetz11Theorem},
	the resolution also works over $\Z$ (i.e.\ the sequence $0 \to \Z \CC_\bullet \to \F_n^\Z(\Sigma_M) \to 0$
	is exact). 
\end{remark}

The second step in our argument is to prove the (local version of)
the tropical Lefschetz hyperplane section theorem for stable intersections.
Even though tropical section theorems are treated in various sources 
(e.g.\ \cite{AB-FilteredGeometricLattices, ARS-LefschetzSectionTheorems}),
it seems that this particular statement has not been covered. 
It is implicit in \cite{Zha-OrlikSolomonAlgebra} however (as explained below).

\begin{lemma} \label{lem:Lefschetzhyperplane}
  Let $H$ denote the standard hyperplane in $\R^N$. Then, for all $p < n$, we have 
	\[
	  \F_p(H \cdot \Sigma_M) = \F_p(\Sigma_M).
	\]
\end{lemma}

\begin{proof}
  It suffices to prove $\F_p(H^{n-p} \cdot \Sigma_M) = \F_p(\Sigma_M)$. 
	By \autoref{lem:hyperplanesection}, $\F_p(H^{n-p} \cdot \Sigma_M) = \F_p(\Sigma_{M^{\leq p}})$.
	By definition, we have $\F_p(\Sigma_{M^{\leq p}}) \subset \F_p(\Sigma_M)$, and it remains
	to show that $\R \CC^{\leq p}_p \to \F_p(\Sigma_M)$ is surjective (as opposed to 
	$\R \CC_p \to \F_p(\Sigma_M)$ which is surjective by definition). 
	
  An indirect, but short proof, works by comparing dimensions. 
	The dimension of $\F_p(\Sigma_{M^{\leq p}})$ can be computed using \autoref{lem:resmaxdim}.
	Conversely, the dimension of $\F_p(\Sigma_M)$ can be computed using the fact that $\F_\bullet(\Sigma_M)$ is isomorphic to the
	Orlik-Solomon algebra $\text{OS}_\bullet(M)$ by \cite{Zha-OrlikSolomonAlgebra}. Expressing both quantities in terms of the Möbius
	function, we find that they agree. 
	
	It is instructive to give a simple direct proof. Let $\FF \in \CC_p$ be a chain of length $p$.
	Let $\gap(\FF)$ be its gap sequence and let $G \supsetneq F$ denote the piece of $\FF$ corresponding to the last non-zero
	entry of $\gap(\FF)$ (in particular, $\rk(G) \geq \rk(F) +2$), and assume $G \neq E$. 
	We denote by $F_1, \dots, F_k$ the the flats of rank $\rk(F) + 1$ such that $G \supsetneq F_i \supsetneq F$. 
	Finally, we denote by $\FF_i \in \CC_p$ the chains obtained from $\FF$ by removing $G$ and inserting $F_i$. 
	The equation for indicator vectors $v_G + (k-1) v_F = \sum v_{F_i}$ implies the equation for volume elements
	\[
	  V_\FF = \sum_{i=1}^k V_{\FF_i}. 
	\]
	But note that the last non-zero entry of the gap sequences $\gap(\FF_i)$ has moved by one position to the front
	(compared to $\gap(\FF)$). 
	By recursion we can express $V_\FF$ in terms of elements $V_{\FF'}$ with gap sequences $\gap(\FF') = (n-p, 0, \dots, 0)$,
	or equivalently, $\FF' \in \CC^{\leq p}_p$. This proves the claim.
\end{proof}

\begin{proof}[\autoref{thm:res}]
  By \autoref{lem:Lefschetzhyperplane} and \autoref{lem:hyperplanesection}, the statement can be reduced to the
	case $p = n$ after intersecting $n-p$ times	with the standard hyperplane. 
	This case was done in \autoref{lem:resmaxdim}. 
\end{proof}

\subsection{Computing the trace side}

We continue by computing $\Tr(\Psi_*, \F_p(\Sigma_M))$ using the resolution
from \autoref{thm:res}.

\begin{lemma} \label{lemTraceChains}
  Let $M$ be a loopless matroid of rank $n+1$ and let 
	$\Psi \colon \Sigma_M \to \Sigma_M$ be a matroidal automorphism. 
	Then for any $p = 0, \dots, n$, we have
	\[
	  (-1)^p \Tr(\Psi_*, \F_p(\Sigma_M)) = \sum_{\substack{F \in \Fix(L(M)) \\ \rk F \leq p}} \mu^\psi(\emptyset, F).
	\]
\end{lemma}

\begin{proof}
  Given a chain $\FF = (F_i)$, we define $\psi(\FF) = (\psi(F_i))$. 
	This induces linear maps on $\R \CC^{\leq p}_i$ for $i=0, \dots,p$ 
	by permuting the generators. 
	By abuse of notation, we  denote these maps by $\Psi_*$. 
	It is obvious that these maps form a morphism 
	of the sequence from \autoref{thm:res}. 
	Hence, using the Hopf trace lemma	\cite[§9, Theorem 2.1]{GD-FixedPointTheory},
	we get
	\[
	  (-1)^p \Tr(\Psi_*, \F_p(\Sigma_M)) = \sum_{i=0}^p (-1)^i \Tr(\Psi_*, \R \CC^{\leq p}_i).
	\]
	Since $\Psi_* \colon \R \CC^{\leq p}_i \to \R \CC^{\leq p}_i$ is given 
	by a permutation of the generators, its trace is equal to the number 
	of fixed generators, i.e.\ the chains $\FF$ with $\psi(\FF) = \FF$. 
	We denote the set of such chains by $\Fix(\CC^{\leq p}_i)$. Obviously,
	it corresponds to the set of chains (with given length/rank) of the
	lattice $\Fix(L(M))$. 
	Note that for any lattice $[\hat{0}, \hat{1}]$, 
	the Möbius function can be expressed in terms
	of chains by 
	\[
	  \mu(\hat{0}, \hat{1}) = \sum_\FF (-1)^{l(\FF)+1}
	\]
	(e.g.\ \cite[Proposition 2.37]{OT-ArrangementsHyperplanes}). 
	Here, the sum runs through 
	all chains $\FF$ of $[\hat{0}, \hat{1}]$ 
	and the length is measured as usual ($\hat{0}, \hat{1}$ are not counted). 
	Applied to our situation, we obtain
	\begin{equation} \begin{split} 
		(-1)^p \Tr(\Psi_*, \F_p(\Sigma_M)) &= \sum_{i=0}^p (-1)^i |\Fix(\CC^{\leq p}_i)| \\
																			 &= \sum_{\substack{F \in \Fix(L(M)) \\ \rk F \leq p}} \sum_{\substack{\FF \in \Fix(\CC) \\ F_1 = F}} (-1)^{l(\FF)} \\
																			 &= \sum_{\substack{F \in \Fix(L(M)) \\ \rk F \leq p}} \mu^\psi(\emptyset, F).
	\end{split} \end{equation}
	Note the change of exponent from $l(\FF)+1$ to $l(\FF)$ due to the fact that we consider $\FF$ as a chain
	of $\Fix(\CC)$, hence $F_1=F$ is counted. 
\end{proof}

It is now easy to finish the proof. 

\begin{proof}[\autoref{eq:secondhalf}]
  By \autoref{lemTraceChains}, we have
	\begin{equation} \begin{split} 
		\sum_{p=0}^n (-1)^{p} \Tr(\Psi_*, \F_p(\Sigma_M)) &= \sum_{p=0}^n \sum_{\substack{F \in \Fix(L(M)) \\ \rk F \leq p}} \mu^\psi(\emptyset, F) \\
		  &= \sum_{F \in \Fix(L(M))} \mu^\psi(\emptyset, F) (n+1 -\rk(F)) \\
			&= - \sum_{F \in \Fix(L(M))} \mu^\psi(\emptyset, F) \rk(F) \\
			&= (-1)^n \beta(\Fix(L(M))).
	\end{split} \end{equation}
	In the second last step, we use (again) that summing the Möbius function over 
	a non-trivial interval gives zero.
\end{proof}

\printbibliography

\subsection*{Contact}

  Johannes Rau, 
	Universidad de los Andes, Carrera 1 No.\ 18A - 12, Bogotá, Colombia; 
  \href{mailto:j.rau@uniandes.edu.co}{j.rau@uniandes.edu.co}.

\subsection*{Data Availability}

Data sharing not applicable to this article as no datasets were generated or analysed during the current study.

\end {document}